\documentclass{commat}

\title{%
    On $\lambda$-pseudo bi-starlike functions related  with Fibonacci numbers
    }

\author{%
    Kaliyappan Vijaya, Gangadharan Murugusundaramoorthy and Hatun \"{O}zlem G\"{u}ney
    }

\affiliation{
    \address{K.~Vijaya --
    School of Advanced Sciences, Vellore Institute of Technology,  Vellore-632014,
India.
        }
    \email{%
    kvijaya@vit.ac.in
    }
    \address{G.~Murugusundaramoorthy --
    School of Advanced Sciences, Vellore Institute of Technology,  Vellore-632014,
India.
        }
    \email{%
    gmsmoorthy@yahoo.com
    }
    \address{H.~G\"{u}ney --
    Dicle University, Faculty of Science, Department of Mathematics, Diyarbak\i r, Turkey
        }
    \email{%
    ozlemg@dicle.edu.tr
    }
    }

\abstract{%
    In this paper we  define a new subclass $\lambda$-bi-pseudo-starlike functions of $\Sigma$  related to
shell-like curves connected with Fibonacci numbers and determine the initial Taylor-Maclaurin
coefficients $|a_2|$ and $|a_3|$ for $f\in\mathcal{PSL}_{\Sigma}^\lambda(\tilde{p}(z)).$ Further we
determine the Fekete-Szeg\"{o} result for the function class
$\mathcal{PSL}_{\Sigma}^\lambda(\tilde{p}(z))$ and for special cases, corollaries are stated which some
of them are new and have not been studied so far.
    }



\VOLUME{32}
\NUMBER{1}
\YEAR{2024}
\firstpage{117}
\DOI{https://doi.org/10.46298/cm.10870}

\begin{paper}

\section{Introduction}

Let $\mathcal{A}$ denote the class of functions $f$  which are {\it
analytic} in the open unit disk
\[
\mathbb{U} =\left\{ z\,
:\,\,z\in\mathbb{C}\, \text{and}\, |z| <1\right\}.
\]
Also let
$\mathcal{S}$ denote the class of functions in $\mathcal{A}$ which
are univalent in $ \mathbb{U}$ and normalised by the conditions
$f(0)=f'(0)-1=0$ and are of the form:
\begin{equation}\label{(1.1)}
    f\left( z\right) =z+\sum_{n=2}^{\infty }a_{n}z^{n}.
\end{equation}
The Koebe one quarter theorem \cite{4} ensures that the image of $ \mathbb{U}$ under
every univalent function $f \in \mathcal{A}$ contains a disk of radius $\frac{1}{4}.$
Thus every univalent function $f$ has an inverse $f^{-1}$ satisfying $$f^{-1}(f(z)) = z,
\ (z \in  \mathbb{U})~{\rm and}~~f(f^{-1}(w)) = w\left(|w| < r_{0}(f), \ r_{0}(f)\geq \frac{1}{4}\right).$$
A function $f \in \mathcal{A}$ is said to be bi-univalent in $ \mathbb{U}$ if both
$f$ and $f^{-1}$ are univalent in $ \mathbb{U}.$ Let $\Sigma$ denote the class of
bi-univalent functions defined in the unit disk $ \mathbb{U}.$ Since $f \in \Sigma$
has the Maclaurian series given by \eqref{(1.1)}, a computation shows that its inverse
$g = f^{-1}$ has the expansion
\begin{equation}\label{(1.2)}
g(w) = f^{-1}(w) = w - a_{2} w^2 + (2a_{2}^2 - a_{3})w^3 + \cdots.
\end{equation}
One can see a short history and examples of functions in the class $\Sigma$ in \cite{14}.
Several authors have  introduced and investigated subclasses of bi-univalent
functions  and obtained bounds for the initial coefficients (see
\cite{2}, \cite{3}, \cite{10}, \cite{14}, \cite{15}, \cite{16}).
For a brief historical account and for several interesting
examples of functions in the class $\Sigma;$ see the pioneering
work on this subject by Srivastava {\it et al.} \cite{14}, which
actually revived the study of bi-univalent functions in
recent years. From the work of Srivastava {\it et al.} \cite{14}, we
choose to recall the following examples of functions in the
class $\Sigma:$
$$\frac{z}{1-z},\quad -\log (1-z)\qquad \text{and} \qquad
\frac{1}{2}\log\left(\frac{1+z}{1-z}\right).$$
We notice that the class $\Sigma$ is not empty. However, the Koebe
function is not a member of~$\Sigma.$

 An analytic function $f$ is said to be subordinate to an analytic function $F$ in $\mathbb{U}$, which we denote by
$f\prec F~( z\in \mathbb{U}) ,$ provided there is an
analytic function $\omega$ defined on $\mathbb{U} $ with $\omega(0)=0$ and $\left\vert
\omega(z)\right\vert <1$ satisfying $f(z)=F(\omega(z)).$ It follows from Schwarz Lemma that $$f(z)\prec
F(z)\quad\Longleftrightarrow f(0)=F(0)~{\rm and}~~ f(\mathbb{U})\subset F(\mathbb{U})~,z\in \mathbb{U}$$
(for details see \cite{4}, \cite{9}).
We recall important subclasses of $\mathcal{S}$ in
geometric function theory such that
if $f\in\mathcal{A}$ and
$$ \frac{zf'(z)}{f(z)}\prec p(z)\quad{\rm and}\quad
1+\frac{zf''(z)}{f'(z)}\prec p(z),
$$ where $p(z)=\frac{1+z}{1-z}$, then we say that $f$ is
starlike and convex, respectively. These
functions form known classes  denoted by $\mathcal{S^{*}}$ and
$\mathcal{C}$, respectively. Recently, in \cite{13}, Sok\'{o}{\l}  introduced the class $\mathcal{SL}$
of shell-like functions as the set of functions $f\in\mathcal{A}$
which is described in the following definition:

\begin{definition} The function $f\in\mathcal{A}$ belongs to the
class $\mathcal{SL}$ if it satisfies the condition that
$$\frac{zf'(z)}{f(z)}\prec\tilde{p}(z)$$ with
$$\tilde{p}(z)=\frac{1+\tau^{2}z^{2}}{1-\tau z-\tau^{2}z^{2}},$$ where $\tau=(1-\sqrt{5})/2\approx -0.618$.
\end{definition}

It should be observed $\mathcal{SL}$ is a subclass of the starlike functions $\mathcal{S^{*}}.$

The function $\tilde{p}$ is not univalent in $\mathbb{U}$, but it is univalent in
the disc of radius $(3-\sqrt{5})/2$. For example, $
\tilde{p}(0)=\tilde{p}(-1/2\tau)=1$ and $\tilde{p}(e^{\mp
i\arccos(1/4)})=\sqrt{5}/5$, and it may also be noticed that
$$\frac{1}{|\tau|}=\frac{|\tau|}{1-|\tau|},$$ which shows that the
number $|\tau|$ divides $[0,1]$ such that it fulfils the golden
section. The image of the unit circle $|z|=1$ under $\tilde{p}$ is a
curve described by the equation given by
$$(10x-\sqrt{5})y^{2}=(\sqrt{5}-2x)(\sqrt{5}x-1)^{2},$$ which is
translated and revolved trisectrix of Maclaurin. The curve
$\tilde{p}(re^{it})$ is a closed curve without any loops for
$0<r\leq r_{0}=(3-\sqrt{5})/2\approx0.38$. For $r_{0}<r<1$, it has a
loop, and for $r=1$, it has a vertical asymptote. Since $\tau$
satisfies the equation $\tau^{2}=1+\tau,$ this expression can be
used to obtain higher powers $\tau^{n}$ as a linear function of
lower powers, which in turn can be decomposed all the way down to a
linear combination of $\tau$ and $1$. The resulting recurrence
relationships yield Fibonacci numbers $u_n$: $$
\tau^{n}=u_{n}\tau+u_{n-1}.$$

In \cite{12} Raina and Sok\'{o}{\l} showed that
\begin{align}
\tilde{p}(z)
&=\frac{1+\tau^{2}z^{2}}{1-\tau z-\tau^{2}z^{2}} \nonumber \\
&=\left(t+\frac{1}{t}\right)\frac{t}{1-t-t^{2}} \nonumber \\
&=\frac{1}{\sqrt{5}} \left(t+\frac{1}{t}\right)\left(\frac{1}{1-(1-\tau)t}-\frac{1}{1-\tau t}\right) \nonumber \\
&=\left(t+\frac{1}{t}\right)\sum_{n=1}^{\infty}u_{n}t^{n} \nonumber \\
&=1+\sum_{n=1}^{\infty}(u_{n-1}+u_{n+1})\tau^{n}z^{n}, \label{(2.6)}
\end{align}
where
\begin{equation}\label{(2.7a)}
u_{n}=\frac{(1-\tau)^{n}-\tau^{n}}{\sqrt{5}} ,
\tau=\frac{1-\sqrt{5}}{2} \,\,\,(n=1,2,\ldots).
\end{equation}
This shows that the relevant connection of $\tilde{p}$ with the
sequence of Fibonacci numbers ${u_{n}}$, such that $u_{0}=0,
\,u_{1}=1,\, u_{n+2}=u_{n}+u_{n+1}$ for $n=0,1,2,\cdots $. And they got

\begin{eqnarray}\label{(2.8)}\tilde{p}(z)&=&1+\sum_{n=1}^{\infty}\tilde{p}_{n}z^{n}\nonumber\\
&=&1+(u_{0}+u_{2})\tau
z+(u_{1}+u_{3})\tau^{2}z^{2}+\sum_{n=3}^{\infty}(u_{n-3}+u_{n-2}+u_{n-1}+u_{n})\tau^{n}z^{n}\nonumber\\
&=&1+\tau z+3\tau^{2}z^{2}+4\tau^{3}z^{3}+7\tau^{4}z^{4}+11\tau^{5}z^{5}+\cdots.
\end{eqnarray}

Let $\mathcal{P}(\beta)$, $0\leq\beta<1$, denote the class of
analytic functions $p$ in $\mathbb{U}$ with $p(0)=1$ and
$Re\{p(z)\}>\beta$. Especially, we will use $\mathcal{P}$ instead of $\mathcal{P}(0).$

\begin{theorem}[\cite{6}] The function
$\tilde{p}(z)=\frac{1+\tau^{2}z^{2}}{1-\tau z-\tau^{2}z^{2}}$
belongs to the class $\mathcal{P}(\sqrt{5}/10)$, where
$\sqrt{5}/10 \approx 0.2236$.
\end{theorem}

Now we recall the following lemma which will be relevant for our study.

\begin{lemma}[\cite{11}]\label{l2}
Let $p\in \mathcal{P}$ with $p(z)=1+c_{1}z+c_{2}z^{2}+\cdots$,
then
\begin{equation}\label{1l2}
    |c_{n}|\leq 2 \qquad for \qquad  n\geq1.
\end{equation}
\end{lemma}

In this present work, we introduce a new subclass of ${\Sigma}$  associated with shell-like functions
connected with Fibonacci numbers and obtain the initial Taylor coefficients $|a_{2}|$ and $|a_{3}|$ for
this function class. Also, we give bounds for the Fekete-Szeg\"{o} functional $|a_{3}-\mu a_{2}^{2}|$ for
this class.

\section{Bi-Univalent function class $\mathcal{PSL}^\lambda_{\Sigma}(\tilde{p}(z))$}

In this section, we introduce a new subclass of ${\Sigma}$  associated with shell-like functions connected
with Fibonacci numbers and obtain the initial Taylor coefficients
$|a_{2}|$ and $|a_{3}|$ for the function class by subordination.

Firstly, let $p(z)=1+p_{1}z+p_{2}z^{2}+\cdots$, and
$p\prec \tilde{p}.$ Then there exists an analytic function $u$ such that $|u(z)|<1$ in $\mathbb{U}$ and
$p(z)=\tilde{p}(u(z))$. Therefore, the function
\begin{equation}\label{1th2.2}
h(z)=\frac{1+u(z)}{1-u(z)}=1+c_{1}z+c_{2}z^{2}+\ldots
\end{equation}
is in the class $\mathcal P(0)$. It follows that
\begin{equation}\label{(2.4)}
u(z)=\frac{c_{1}z}{2}+\left(c_{2}-\frac{c_{1}^{2}}{2}\right)\frac{z^{2}}{2}+\left(c_{3}-c_{1}c_{2}+\frac{c_{1}^{3}}{4}\right)\frac{z^{3}}{2}+\cdots
\end{equation}
and
\begin{eqnarray}\label{(2.7)}
\tilde{p}(u(z))&=&
1+\frac{\tilde{p}_{1}c_{1}z}{2}+
\left\{\frac{1}{2}\left(c_{2}-\frac{c_{1}^{2}}{2}\right)\tilde{p}_{1}+\frac{c_{1}^{2}}{4}\tilde{p}_{2}\right\}z^{2}\nonumber\\
&~&+
 \left\{\frac{1}{2}\left(c_{3}-c_{1}c_{2}+\frac{c_{1}^{3}}{4}\right)\tilde{p}_{1}
 +\frac{1}{2}c_{1}\left(c_{2}-\frac{c_{1}^{2}}{2}\right)\tilde{p}_{2}+\frac{c_{1}^{3}}{8}\tilde{p}_{3}\right\}z^{3}+\cdots.
\end{eqnarray}
And similarly,
there exists an analytic function $v$ such that $|v(w)|<1$ in $\mathbb{U}$ and such that
$p(w)=\tilde{p}(v(w))$. Therefore, the function
\begin{equation}\label{1th2.2v}
k(w)=\frac{1+v(w)}{1-v(w)}=1+d_{1}w+d_{2}w^{2}+\ldots
\end{equation}
is in the class $\mathcal P(0)$. It follows that
\begin{equation}\label{(2.10)}
v(w)=\frac{d_{1}w}{2}+\left(d_{2}-\frac{d_{1}^{2}}{2}\right)\frac{w^{2}}{2}+\left(d_{3}-d_{1}d_{2}+\frac{d_{1}^{3}}{4}\right)\frac{w^{3}}{2}+\cdots
\end{equation}
and
\begin{eqnarray}\label{(2.5)}
\tilde{p}(v(w))&=&1+\frac{\tilde{p}_{1}d_{1}w}{2}+
\left\{\frac{1}{2}\left(d_{2}-\frac{d_{1}^{2}}{2}\right)\tilde{p}_{1}
+\frac{d_{1}^{2}}{4}\tilde{p}_{2}\right\}w^{2}\nonumber\\
&~&+
 \left\{\frac{1}{2}\left(d_{3}-d_{1}d_{2}+\frac{d_{1}^{3}}{4}\right)\tilde{p}_{1}
 +\frac{1}{2}d_{1}\left(d_{2}-\frac{d_{1}^{2}}{2}\right)\tilde{p}_{2}
 +\frac{d_{1}^{3}}{8}\tilde{p}_{3}\right\}w^{3}+\cdots.
\end{eqnarray}
  The class $\mathcal{L}_\lambda(\alpha)$ of $\lambda$-pseudo-starlike functions
of order $\alpha$ was introduced and investigated
by Babalola \cite{1}.  A function $f\in\mathcal A$ is in the
class $\mathcal{L}_\lambda(\alpha)$ if it satisfies
\begin{equation*}\label{1.4}
    \Re\left(\frac{z(f'(z))^\lambda}{f(z)}\right) > \alpha,\quad (0\leq \alpha<1)
\end{equation*} where $\lambda \geq 1,\lambda \in \mathbb{R}$ and $z \in \mathbb{U}$.
In \cite{1} it was showed that all pseudo-starlike functions are
Bazilevi\v{c} functions of type $\left(1-1/\lambda\right)$ and of
order $\alpha^{1/\lambda}$ and univalent in open unit disk
$\mathbb{U}$.

 Recently Joshi et al. \cite{5}  defined the bi-pseudo-starlike functions class and
obtained the bounds for the initial coefficients $|a_2|$ and
$|a_3|$. In this paper we  define a new class
$\mathcal{PSL}^\lambda_{\Sigma}(\tilde{p}(z)),$$\lambda$-bi-pseudo-starlike
functions of $\Sigma$ related to shell-like curves connected with Fibonacci numbers and determine the
bounds for the initial Taylor-Maclaurin coefficients of $|a_2|$ and $|a_3|.$  Further we consider the
Fekete-Szeg\"{o} problem in this class and the special cases are stated as corollaries which are new and
have not been studied so far.
\begin{definition}\label{def1} For $\lambda\geq 1$ and $\lambda\in \mathbb{R},$ a  function $f \in \Sigma$
of the
form (\ref{(1.1)}) is said to be in the class $\mathcal{PSL}^\lambda_{\Sigma}(\tilde{p}(z))$
if the following subordination hold:
\begin{equation}\label{(1.13)}
\frac{z(f'(z))^\lambda}{f(z)}\prec\tilde{p}(z)=\frac{1+\tau^{2}z^{2}}{1-\tau z-\tau^{2}z^{2}}
\end{equation}
and
\begin{equation}\label{(1.14)}
\frac{w(g'(w))^\lambda}{g(w)}\prec\tilde{p}(w)=\frac{1+\tau^{2}w^{2}}{1-\tau w-\tau^{2}w^{2}}
\end{equation}
where $\tau=(1-\sqrt{5})/2\approx -0.618$ where $z,w \in \mathbb{U} $ and $g$ is given by (\ref{(1.2)}).
\end{definition}
Specialising the parameter $\lambda=1$ and $\lambda=2,$ we  have the following remarks, respectively:
\begin{remark}[\cite{8}] \label{def2}
For $\lambda=1$ a function $f\in \Sigma$
is  in the class $\mathcal{PSL}^1_{\Sigma}(\tilde{p}(z)) \equiv \mathcal{SL}_{\Sigma}(\tilde{p}(z))$
if
the following conditions are satisfied:
\begin{equation}\label{(1.9)}
\frac{zf'(z)}{f(z)}\prec\tilde{p}(z)=\frac{1+\tau^{2}z^{2}}{1-\tau z-\tau^{2}z^{2}}
\end{equation}
and
\begin{equation}\label{(1.10)}
\frac{wg'(w)}{g(w)}\prec\tilde{p}(w)=\frac{1+\tau^{2}w^{2}}{1-\tau w-\tau^{2}w^{2}},
\end{equation}
where $\tau=(1-\sqrt{5})/2\approx -0.618$ where $z,w \in \mathbb{U} $ and $g$ is given by (\ref{(1.2)}).
\end{remark}
\begin{remark}
For $\lambda=2$ a function $f\in \Sigma$
is  in the class $\mathcal{PSL}^2_{\Sigma}(\tilde{p}(z)) \equiv \mathcal{GSL}_{\Sigma}(\tilde{p}(z))$
if
the following conditions are satisfied:
\begin{equation}\label{(1.11)}
\left(f'(z)\frac{zf'(z)}{f(z)}\right)\prec\tilde{p}(z)=\frac{1+\tau^{2}z^{2}}{1-\tau z-\tau^{2}z^{2}}
\end{equation}
and
\begin{equation}\label{(1.12)}
\left(g'(w)\frac{wg'(w)}{g(w)}\right)\prec\tilde{p}(w)=\frac{1+\tau^{2}w^{2}}{1-\tau w-\tau^{2}w^{2}},
\end{equation}
where $\tau=(1-\sqrt{5})/2\approx -0.618$ where $z,w \in \mathbb{U} $ and $g$ is given by (\ref{(1.2)}).
\end{remark}

In the following theorem we determine the initial Taylor coefficients $|a_{2}|$ and $|a_{3}|$ for the
function class $\mathcal{PSL}^\lambda_{\Sigma}(\tilde{p}(z)).$ Later we will reduce these bounds to other
classes for special cases.
\begin{theorem}\label{thm1}
Let $f$ given by \eqref{(1.1)} be in the class $\mathcal{PSL}^\lambda_{\Sigma}(\tilde{p}(z)),$ then
\begin{equation}\label{c7e3.3}
|a_{2}|\leq\frac{|\tau|}{\sqrt{(2\lambda-1)^2-(10\lambda^2-11\lambda+3)\tau}}
\end{equation}
and
\begin{equation}\label{c7e3.4}
|a_{3}|\leq\frac{|\tau|\left[(2\lambda-1)^{2}-2(5\lambda^{2}-4\lambda+1)\tau\right]}
{(3\lambda-1)\left[(2\lambda-1)^{2}-(10\lambda^{2}-11\lambda+3)\tau\right]},
\end{equation}where $\lambda \geq 1 .$
\end{theorem}
\begin{proof}
Let $f \in \mathcal{PSL}^\lambda_{\Sigma}(\tilde{p}(z))$ and $g = f^{-1}.$ Considering \eqref{(1.13)} and
\eqref{(1.14)}, we have
\begin{equation}
\frac{z(f'(z))^\lambda}{f(z)}=\tilde{p}(u(z))
\end{equation}
and
\begin{equation}
\frac{w(g'(w))^\lambda}{g(w)}=\tilde{p}(v(w)),
\end{equation}
where $\tau=(1-\sqrt{5})/2\approx -0.618$ where $z,w \in \mathbb{U} $ and $g$ is given by (\ref{(1.2)}).
Since
\[\frac{z(f'(z))^\lambda}{f(z)}=1+(2\lambda-1)a_2z+[(3\lambda-1)a_3+\left(2\lambda^2-4\lambda+1\right)a_2^2]z^2
+\dots\]
 and \[\frac{w(g'(w))^\lambda}{g(w)}=
 1-(2\lambda-1)a_2w+[\left(2\lambda^2+2\lambda-1\right)a_2^2-(3\lambda-1)a_3]w^2+\dots.\]
Thus we have
\begin{align}\label{2.22}
1+(2\lambda-1)a_2z&+[(3\lambda-1)a_3+\left(2\lambda^2-4\lambda+1\right)a_2^2]z^2
+\dots \\ ={ }&{ }1+\frac{\tilde{p}_{1}c_{1}}{2}z+
\left[\frac{1}{2}\left(c_{2}-\frac{c_{1}^{2}}{2}\right)\tilde{p}_{1}
+\frac{c_{1}^{2}}{4}\tilde{p}_{2}\right]z^{2}\nonumber\\
&+
 \left[\frac{1}{2}\left(c_{3}-c_{1}c_{2}+\frac{c_{1}^{3}}{4}\right)\tilde{p}_{1}
 +\frac{1}{2}c_{1}\left(c_{2}-\frac{c_{1}^{2}}{2}\right)\tilde{p}_{2}
 +\frac{c_{1}^{3}}{8}\tilde{p}_{3}\right]z^{3}+ \cdots.
\end{align}
and
\begin{align}\label{2.22-1}\nonumber
1-(2\lambda-1)a_2w&+[\left(2\lambda^2+2\lambda-1\right)a_2^2-(3\lambda-1)a_3]w^2
+\dotsb ,\\
={ }&{ }1+\frac{\tilde{p}_{1}d_{1}w}{2}+
\left[\frac{1}{2}\left(d_{2}-\frac{d_{1}^{2}}{2}\right)\tilde{p}_{1}+\frac{d_{1}^{2}}{4}\tilde{p}_{2}\right]w^{2}\nonumber\\
&+
 \left[\frac{1}{2}\left(d_{3}-d_{1}d_{2}+\frac{d_{1}^{3}}{4}\right)\tilde{p}_{1}
 +\frac{1}{2}d_{1}\left(d_{2}-\frac{d_{1}^{2}}{2}\right)\tilde{p}_{2}+\frac{d_{1}^{3}}{8}\tilde{p}_{3}\right]w^{3}+\cdots.
\end{align}
It follows from (\ref{2.22}) and (\ref{2.22-1}) that
\begin{equation}\label{2.24}
(2\lambda-1)a_2=\frac{c_{1}\tau}{2},
\end{equation}
\begin{equation}\label{2.25}
(3\lambda-1)a_3+\left(2\lambda^2-4\lambda+1\right)a_2^2=\frac{1}{2}\left(c_{2}-\frac{c_{1}^{2}}{2}\right)\tau+\frac{c_{1}^{2}}{4}3\tau^{2},
\end{equation}
and
\begin{equation}\label{2.26}
-(2\lambda-1)a_2=\frac{d_{1}\tau}{2},
\end{equation}
\begin{equation}\label{2.27}
\left(2\lambda^2+2\lambda-1\right)a_2^2-(3\lambda-1)a_3=\frac{1}{2}\left(d_{2}-\frac{d_{1}^{2}}{2}\right)\tau+\frac{d_{1}^{2}}{4}3\tau^{2}.
\end{equation}
From (\ref{2.24}) and (\ref{2.26}), we have
\begin{equation}\label{2.28}
c_{1}=-d_{1},
\end{equation}
and
\begin{equation}\label{2.29}
a_{2}^{2}=\frac{(c_{1}^{2}+d_{1}^{2})}{8(2\lambda-1)^2}\tau^{2}.
\end{equation}
Hence
\begin{equation}\label{2.29a}
|a_{2}|\leq\frac{|\tau|}{2\lambda-1}.
\end{equation}
Now, by summing (\ref{2.25}) and (\ref{2.27}), we obtain
\begin{equation}\label{2.30}
2(2\lambda-1)^2a_{2}^{2}=\frac{1}{2}(c_{2}+d_{2})\tau-\frac{1}{4}(c_{1}^{2}+d_{1}^{2})\tau+\frac{3}{4}(c_{1}^{2}+d_{1}^{2})\tau^{2}.
\end{equation}
Substituting  (\ref{2.29}) in (\ref{2.30}), we have
\begin{equation}\label{(2.35)}
2\left[(2\lambda-1)^2-(10\lambda^2-11\lambda+3)\tau\right]a_{2}^{2}=\frac{1}{2}(c_{2}+d_{2})\tau^2.
\end{equation}
Therefore, using Lemma (\ref{l2}) we obtain
\begin{equation}\label{2.32}
|a_{2}| \leq \frac {|\tau| } { \sqrt{(2\lambda-1)^2-(10\lambda^2-11\lambda+3)\tau}}.
\end{equation}
Now, so as to find the bound on $|a_{3}|,$ let's subtract from (\ref{2.25}) and (\ref{2.27}). So, we find
\begin{equation}\label{(2.37)}
2(3\lambda-1)a_{3}-2(3\lambda-1)a_{2}^{2}=\frac{1}{2}\left(c_{2}-d_{2}\right)\tau.
\end{equation}
Hence, we get
\begin{equation}\label{2.33}
2(3\lambda-1)|a_{3}|\leq 2|\tau|+ 2(3\lambda-1) |a_{2}|^{2}.
\end{equation}
Then, in view of  (\ref{2.32}), we obtain
\begin{equation*}
|a_{3}|\leq\frac{|\tau|\left[(2\lambda-1)^{2}-2(5\lambda^{2}-4\lambda+1)\tau\right]}{(3\lambda-1)\left[(2\lambda-1)^{2}-(10\lambda^{2}-11\lambda+3)\tau\right]}.
\qedhere
\end{equation*}
\end{proof}

By taking the parameter $\lambda =1$ and  $\lambda =2$ in the above theorem, we  have the following the
initial Taylor coefficients $|a_{2}|$ and $|a_{3}|$ for the function classes
$\mathcal{SL}_{\Sigma}(\tilde{p}(z))$ and $\mathcal{GSL}_{\Sigma}(\tilde{p}(z)),$ respectively.

\begin{corollary}[\cite{8}] \label{2.2}
Let $f$ given by \eqref{(1.1)} be in the class $\mathcal{SL}_{\Sigma}(\tilde{p}(z)),$ then
\begin{equation}\label{c7e3.3a}
|a_{2}| \leq \frac {|\tau| } { \sqrt{1-2\tau}}
\end{equation}
and
\begin{equation}\label{c7e3.4a}
|a_{3}|\leq\frac {|\tau|(1-4\tau)} {2(1-2\tau)}.
\end{equation}
\end{corollary}

\begin{corollary}\label{2.3}
Let $f$ given by \eqref{(1.1)} be in the class $\mathcal{GSL}_{\Sigma}(\tilde{p}(z)),$ then
\begin{equation}\label{c7e3.3b}
|a_{2}|\leq\frac{|\tau|}{\sqrt{9-21\tau}}
\end{equation}
and
\begin{equation}\label{c7e3.4b}
|a_{3}|\leq\frac{|\tau|(9-26\tau)}{5(9-21\tau)}.
\end{equation}
\end{corollary}

\section{Fekete-Szeg\"{o} inequality for the function class $\mathcal{PSL}^\lambda_{\Sigma}(\tilde{p}(z))$}

Fekete and Szeg\"{o} \cite{7} introduced the generalised functional $|a_{3}-\mu a_{2}^{2}|,$
where $\mu$ is some real number. Due to Zaprawa \cite {17},in the following theorem we determine the
Fekete-Szeg\"{o} functional for $f\in \mathcal{PSL}^\lambda_{\Sigma}(\tilde{p}(z))$.

\begin{theorem}\label{thm2} Let $f$ given by (\ref{(1.1)}) be in the class
$\mathcal{PSL}^\lambda_{\Sigma}(\tilde{p}(z))$ and $\mu\in\mathbb{R},$ then
\begin{equation*}
|a_{3}-\mu a_{2}^{2}|\leq\left\{ \begin{array}{ll}
\frac{|\tau|}{4(3\lambda-1)} , & \mbox{ $0\leq |h(\mu)| \leq \frac{|\tau|}{4(3\lambda-1)}$},\\
4|h(\mu)| , &\,\,\,\,\,\,\,\,\,\,\,\mbox{ $ |h(\mu)|\geq \frac{|\tau|}{4(3\lambda-1)}$},\end{array} \right.
\end{equation*}
where
\begin{equation}\label{(3.1)}
h(\mu)=\frac{(1-\mu)\tau^{2}}{4\left[(2\lambda-1)^2-(10\lambda^2-11\lambda+3)\tau\right]}.
\end{equation}
\end{theorem}
\begin{proof}
From (\ref{(2.35)}) and (\ref{(2.37)}) we obtain
\begin{align}
a_{3}-\mu a_{2}^{2}
={ }&{ }(1-\mu)\frac{\tau^{2}(c_{2}+d_{2})}{4\left[(2\lambda-1)^2-(10\lambda^2-11\lambda+3)\tau\right]}
+\frac{\tau(c_{2}-d_{2})}{4(3\lambda-1)} \label{(3.1c)} \\
={ }&{ }\left(\frac{(1-\mu)\tau^{2}}{4\left[(2\lambda-1)^2-(10\lambda^2-11\lambda+3)\tau\right]}+
\frac{\tau}{4(3\lambda-1)}\right)c_{2} \nonumber \\
&+\left(\frac{(1-\mu)\tau^{2}}{4\left[(2\lambda-1)^2-(10\lambda^2-11\lambda+3)\tau\right]}
-\frac{\tau}{4(3\lambda-1)}\right)d_{2}. \nonumber
\end{align}
So we have
\begin{equation}\label{(3.45)}
a_{3}-\mu a_{2}^{2}=\left(h(\mu)+\frac{|\tau|}{4(3\lambda-1)}\right)c_{2}
+\left(h(\mu)-\frac{|\tau|}{4(3\lambda-1)}\right)d_{2},
\end{equation}
where
\begin{equation}\label{(3.1e)}
h(\mu)=\frac{(1-\mu)\tau^{2}}{4\left[(2\lambda-1)^2-(10\lambda^2-11\lambda+3)\tau\right]}.
\end{equation}

Then, by taking modulus of (\ref{(3.45)}), we conclude that
\[|a_{3}-\mu a_{2}^{2}|\leq\left\{ \begin{array}{ll}
\frac{|\tau|}{4(3\lambda-1)} , & \mbox{ $0\leq |h(\mu)| \leq \frac{|\tau|}{4(3\lambda-1)}$},\\
4|h(\mu)| , &\,\,\,\,\,\,\,\,\,\,\,\mbox{ $ |h(\mu)|\geq \frac{|\tau|}{4(3\lambda-1)}.$}\end{array} \right.
\qedhere
\]
\end{proof}

Taking $\mu=1$, we have the following corollary.
\begin{corollary}\label{c3} If $f\in \mathcal{PSL}^\lambda_{\Sigma}(\tilde{p}(z))$, then
\begin{equation}\label{(3.1aa)}
|a_{3}-a_{2}^{2}|\leq\frac{|\tau|}{4(3\lambda-1)}.
\end{equation}
\end{corollary}

By specialising the parameter $\lambda =1$ and  $\lambda=2$ in the above theorem, we  have the following the Fekete-Szeg\"{o} inequalities for the function classes $\mathcal{SL}_{\Sigma}(\tilde{p}(z))$
and $\mathcal{GSL}_{\Sigma}(\tilde{p}(z)),$ respectively.

\begin{corollary}[\cite{8}] \label{cor2}
Let $f$ given by (\ref{(1.1)}) be in the class
$\mathcal{SL}_{\Sigma}(\tilde{p}(z))$ and $\mu\in\mathbb{R},$. then we have
\begin{equation*}
|a_{3}-\mu a_{2}^{2}|\leq\left\{ \begin{array}{ll}
\frac{|\tau|}{8} , & \mbox{ $0\leq |h(\mu)| \leq \frac{|\tau|}{8}$},\\
4|h(\mu)| , &\,\,\,\,\,\,\,\,\,\,\,\mbox{ $ |h(\mu)|\geq \frac{|\tau|}{8}$},\end{array} \right.
\end{equation*}
where
\begin{equation}\label{(3.1b)}
h(\mu)=\frac{(1-\mu)\tau^{2}}{4\left[1-2\tau\right]}.
\end{equation}
Further if $\mu=1$ we get
$$|a_{3}-a_{2}^{2}|\leq\frac{|\tau|}{8}.$$
\end{corollary}
\begin{corollary}\label{cor3} Let $f$ given by (\ref{(1.1)}) be in the class
$\mathcal{GSL}_{\Sigma}(\tilde{p}(z))$ and $\mu\in\mathbb{R},$ then we have
\begin{equation*}
|a_{3}-\mu a_{2}^{2}|\leq\left\{ \begin{array}{ll}
\frac{|\tau|}{20} , & \mbox{ $0\leq |h(\mu)| \leq \frac{|\tau|}{20}$},\\
4|h(\mu)| , &\,\,\,\,\,\,\,\,\,\,\,\mbox{ $ |h(\mu)|\geq \frac{|\tau|}{20}$},\end{array} \right.
\end{equation*}
where
\begin{equation}\label{(3.1 a)}
h(\mu)=\frac{(1-\mu)\tau^{2}}{4\left[9-21\tau\right]}.
\end{equation}
Further if $\mu=1$ we get
$$|a_{3}-a_{2}^{2}|\leq\frac{|\tau|}{20}.$$
\end{corollary}

\subsection*{Acknowledgements}The authors thank the referees of this paper for their insightful suggestions and
corrections to improve the paper in present form.


\EditInfo{March 04, 2019}{November 09, 2021}{Karl Dilcher}

\end{paper}